\documentclass{elsart3-1}

\usepackage{amssymb, amsmath}

\usepackage[english, french]{babel}

\newtheorem{theorem}{Theorem}[section]
\newtheorem{lemma}[theorem]{Lemma}
\newtheorem{e-proposition}[theorem]{Proposition}

\newtheorem{e-definition}[theorem]{Definition\rm}
\newtheorem{remark}{\it Remark\/}

\renewcommand{\theequation}{\arabic{equation}}
\def\og{\leavevmode\raise.3ex\hbox{$\scriptscriptstyle\langle\!\langle$~}}
\def\fg{\leavevmode\raise.3ex\hbox{~$\!\scriptscriptstyle\,\rangle\!\rangle$}}

\newcommand{\Rey}{\mathrm{Re}}
\newcommand{\Hronde}{{\mathcal H}}
\newcommand{\pd}[2]{\frac{\partial #1}{\partial\/ #2}}

\begin{document}
\centerline{}
\begin{frontmatter}



\selectlanguage{english}
\title{Derivation of an intermediate viscous Serre--Green--Naghdi equation}


\selectlanguage{english}
\author[DDS]{Denys Dutykh}
\ead{Denys.Dutykh@univ-smb.fr}
\author[HLM]{Herv\'e Le Meur},
\ead{Herve.Le.Meur@u-picardie.fr}

\address[DDS]{Univ. Grenoble Alpes, Univ. Savoie Mont Blanc, CNRS, LAMA, 73000 Chamb\'ery, France}
\address[HLM]{LAMFA, CNRS, UMR 7352, Universit\'e Picardie Jules Verne, 80039 Amiens, France}


\medskip
\begin{center}
{\small Received *****; accepted after revision +++++\\
Presented by XXXXXX}
\end{center}

\begin{abstract}
\selectlanguage{english}
In this note we present the current status of the derivation of a viscous Serre--Green--Naghdi system. For this goal, the flow domain is separated into two regions. The upper region is governed by inviscid Euler equations, while the bottom region (the so-called boundary layer) is described by Navier-Stokes equations. We consider a particular regime linking the Reynolds number and the shallowness parameter. The computations presented in this note are performed in the fully nonlinear regime. The boundary layer flow reduces to a Prantdl-like equation. Further approximations seem to be needed to obtain a tractable model.


\vskip 0.5\baselineskip

\selectlanguage{french}
\noindent{\bf R\'esum\'e} \vskip 0.5\baselineskip \noindent
{\bf Obtention des \'equations interm\'ediaires de Serre--Green--Naghdi visqueuse.}
Dans cette note nous pr\'esentons l'\'etat actuel de l'obtention du syst\`eme de type Serre--Green--Naghdi visqueux dans un canal. Nous s\'eparons le domaine fluide en deux couches. La couche sup\'erieure est d\'ecrite par les \'equations d'Euler tandis que la couche limite en-dessous, ob\'eit aux \'equations de Navier-Stokes. Nous consid\'erons un r\'egime pleinement non lin\'eaire o\`u le nombre de Reynolds est li\'e au rapport de la longueur d'onde typique \`a la profondeur moyenne. La dynamique de l'\'ecoulement dans la couche limite se ram\`ene \`a une \'equation de type Prantdl. Des hypoth\`eses suppl\'ementaires sont n\'ecessaires afin d'obtenir un mod\`ele utilisable en pratique.


\end{abstract}
\end{frontmatter}

\selectlanguage{french}
\section*{Version fran\c{c}aise abr\'eg\'ee}
Nous tentons d'obtenir un mod\`ele r\'eduit aux \'equations de
l'\'ecoulement d'un fluide visqueux dans un canal peu profond. Nous ne
supposons pas l'amplitude des vagues comme petite (r\'egime non
lin\'eaire). Comme dans le cas du r\'egime lin\'eaire
(cf. \cite{LeMeur2015}), nous r\'esolvons les \'equations dans la zone
principale gouvern\'ee par des \'equations d'Euler pour arriver \`a
(\ref{eq:15}). Puis nous tentons la m\^eme r\'esolution dans la couche
limite, mais ne pouvons aller plus loin que (\ref{eq:21}). Cette
derni\`ere \'equation est de type Prandtl. Elle est connue pour un
comportement tr\`es sensible \`a chacun de ses termes
(cf. \cite{Gerard-Varet_Dormy_2010}) et donc
laisse peu d'espoir pour \^etre simplifi\'ee afin d'obtenir un
mod\`ele 1D.

\selectlanguage{english}

\section{Introduction}

The water wave theory has been essentially developed in the framework
of the inviscid, and very often also irrotational, Euler
equations. However, various viscous effects are inevitably present in
laboratory experiments and even more in the real world. Thus, the
conservative conventional models have to be supplemented with
dissipative effects to improve the quality of their predictions. A
straightforward energy balance asymptotic analysis shows that the main
dissipation takes place at the bottom boundary layer
\cite[Section~\textsection2]{Dutykh2008a} (or at the lateral walls if
they are also present \cite{Kit1989}). In this way, the corresponding
long wave Boussinesq-type systems have been derived taking into
account the boundary layer effects \cite{Dutykh2014f}. In
\cite{LeMeur2015}, the author derives the viscous Boussinesq model
without the irrotationality assumption. Other articles already took
the vorticity into account, even for fully nonlinear
Boussinesq equations (here called Serre-Green-Naghdi or SGN)
\cite{Lannes_Castro_2014}. Fully nonlinear models
are becoming very popular. In the present note we report the current
status of the derivation of a viscous counterpart of the well-known
SGN equations. The asymptotic regime relates the Reynolds number to
the shallowness parameter.

\section{Primary equations}

Consider the flow of an incompressible liquid in a physical two-dimensional space over a flat bottom and with a free surface. We assume additionally that the fluid is homogeneous (i.e.~the density $\rho\ =\ \mathrm{const}$) and the gravity acceleration $g$ is constant. For the sake of simplicity, in this study we neglect all other forces (such as the Coriolis force and friction). Hence, we deal with pure gravity waves. We introduce a Cartesian coordinate system $O\, \tilde{x}\,\tilde{y}\,$. The horizontal line $O\, \tilde{x}$ coincides with the still water level $\tilde{y} = 0$ and the axis $O\, \tilde{y}$ points vertically upwards. The fluid layer is bounded below by the horizontal solid bottom $\tilde{y} = -d$ and above by the free surface $\tilde{y} = \tilde{\eta}\,(\tilde{x},\,\tilde{t})\,$.

In order to make the equations dimensionless, we choose a characteristic horizontal length $\ell$, vertical height of the free surface $A$ and mean depth $d$. All this enables us to define a characteristic velocity $c_0=\sqrt{g d}$. Then one may define dimensionless independent variables:
\begin{equation*}
  \tilde{x}=\ell x, \; \tilde{y}=d y, \; \tilde{t} = t \ell /c_0.
\end{equation*}
This enables us to define the dimensionless fields:
\begin{equation*}
  \tilde{u}=c_0 u, \; \tilde{v}=\frac{d c_0}{\ell} v, \; \tilde{p}=\tilde{p}_{\mathrm{atm}} - \rho g d \, y+\rho g d \, p, \; \tilde{\eta}(\tilde{x},\tilde{y},\tilde{t})= A\eta(x,y,t).
\end{equation*}
We also define some dimensionless numbers, characteristic of the flow:
\begin{equation*}
  \varepsilon =\frac{A}{d}, \; \mu^2=\frac{d^2}{\ell^2}, \; \Rey=\frac{\rho c_0 d}{\nu}.
\end{equation*}
The system of Navier-Stokes equations can then be written in 2D and in dimensionless variables:
\begin{equation}
  \label{eq:NS_adim}
  \left\{
  \begin{array}{rcl}
    u_t+u u_x+vu_y- 1/\Rey\left(\mu u_{xx}+ u_{yy}/\mu \right)+p_x & = 0 & \\
    \mu^2(v_t+u v_x+v v_y)-\mu^2/\Rey\left(\mu v_{xx}+ v_{yy}/\mu\right)+p_y & = 0 & \\
    u_x+v_y & = 0 & \\
    \left[-\left(p- \varepsilon \eta \right) {\boldsymbol I}+\displaystyle \frac{2}{\Rey}\left( \begin{array}{cc} \mu u_x & (u_y+\mu^2v_x)/2 \\(u_y+\mu^2v_x)/2  & \mu v_y \end{array} \right) \right]{\Bigg |}_{\varepsilon \eta} {\boldsymbol n} & = 0 & \mbox{ on }y=\varepsilon \eta\\
    \eta_t + u(y=\varepsilon \eta) \eta_x -v(y=\varepsilon \eta)/\varepsilon& =0  & \mbox{ on }y=\varepsilon \eta\\
    u(y=-1)=v(y=-1) & = 0 & ,
  \end{array}
  \right.
\end{equation}
where we denote $u|_{\varepsilon \eta}=u(y=\varepsilon \eta)=u(x,y=\varepsilon \eta(x,t),t)$.

One could assume the fields to be small around the hydrostatic flow (which is lifted by the change of field from $\tilde{p}$ to $p$), so around $(u,v,p,\eta) \simeq 0$. Such an assumption is contradictory with our nonlinear assumption where $\varepsilon$ is assumed not to be small. Yet, should we make this assumption, we would be led to a linear system identical (up to changes of variables) to System~(7) of \cite{LeMeur2015}. The study of the linear regime suggests to assume, not only in the Boussinesq regime:
\begin{equation}\label{assumption_Re}
  \Rey \simeq \mu^{-5}\,.
\end{equation}

Below, we solve the problem in the bulk part where Euler's equations are justified to apply (Section \ref{sec.UP}), then try to solve the velocity in the boundary layer (Section \ref{sec.BL}). In this last section, we are led to Prandtl's equation that prohibits any further advance to the best of our knowledge.

[??? \`a nettoyer] What is the size of the boundary layer where the no-slip condition forces the fluid to have a large gradient of velocity ? In the same way as in \cite{LeMeur2015}, one may assume it is of size $\mu^2$:
\begin{equation}\label{eq6.55}
  y=-1+\mu^2 \, \gamma.
\end{equation}

One might be surprised that the gravity-viscosity layer be of size $O(\mu^2)$ (or a little larger) while one usually assumes the size of the viscous layer to be of size $O(\Rey^{-1/2})=O(\mu^{5/2})$. Indeed the classical term stems from the $1/\Rey \, u_{yy}$ term which is replaced here by $1/(\mu \Rey) \, u_{yy}$. So $\mu \Rey_{us}=\Rey_{classical}$ and the size of the boundary layer is $\Rey_{classical}^{-1/2}=(\mu\Rey_{us})^{-1/2}=(\mu^{-4})^{-1/2}=\mu^2$!

\subsection{Resolution in the upper part (Euler)}
\label{sec.UP}

In the upper part, $y \gg -1+ \mu^2$ and $\mu^4$ is small. So one may drop the Laplacian and keep from (\ref{eq:NS_adim}):
\begin{equation}\label{eq:8}
\left\{
\begin{aligned}
  u_t+u u_x +v u_y+p_x & = O\left(u_{yy}/(\mu \Rey)\right)+O(\mu^6) & &\\
  \mu^2(v_t+u v_x +v v_y)+p_y & = O(\mu^6) & &\\
  u_x+v_y & =0 & \\
  -p+\varepsilon \eta & = O(u_y|_{\varepsilon \eta}/\Rey )+O(\mu/\Rey) & &\mbox{ on } y=\varepsilon \eta\\
  (p-\varepsilon\eta)\varepsilon \eta_x+2(-\mu u_x\varepsilon\eta_x+(u_y+\mu^2 v_x))/\Rey & = 0 & &\mbox{ on } y=\varepsilon \eta\\
  \varepsilon\left(\eta_t + u|_{\varepsilon\eta} \eta_x \right) & =v|_{\varepsilon\eta}& & \mbox{ on } y=\varepsilon \eta.
\end{aligned}
\right.
\end{equation}
First, one may notice that the viscosity terms are no more present inside the domain. It is argued in \cite{LeMeur2015} that one may (and even must) then drop the fifth equation from this system, due to the fact that the fluid is indeed no more viscous in this part of the domain.

It is classical to use (\ref{eq:8})$_3$ to get
\begin{equation}
  \label{eq:9}
  v=v|_{\varepsilon\eta}-\int_{\varepsilon \eta}^y u_x \, {\rm d} y',
\end{equation}
where $v|_{y=\varepsilon\eta}$ is given by (\ref{eq:8})$_6$.
One may use this vertical velocity in (\ref{eq:8})$_2$ to compute $p_y$.
Thanks to (\ref{eq:8})$_4$, one has:
\begin{align}
  \nonumber
p= &\varepsilon \eta -\mu^2 \Bigg[(y-\varepsilon \eta)\left( (v|_{\varepsilon\eta})_t+u_x|_{\varepsilon\eta} \varepsilon \eta_t \right)+\left(\int_{\varepsilon \eta}^y u \right) \left( (v|_{\varepsilon\eta})_x+u_x|_{\varepsilon\eta} \varepsilon \eta_x -\left(\int_{\varepsilon \eta}^y u_x \right) v|_{\varepsilon \eta}\right)\Bigg.\\
\label{eq:10} & \left. -\int_{\varepsilon \eta}^y\int_{\varepsilon \eta}^{y'}u_{xt}-\int_{\varepsilon \eta}^y\left( u\int_{\varepsilon \eta}^{y'}u_{xx} \right) +\int_{\varepsilon \eta}^y\left(u_x \int_{\varepsilon \eta}^{y'} u_x \right)\right]+O\left(\frac{u_y|_{\varepsilon \eta}}{\Rey} \right)+O\left( \frac{\mu}{\Rey} \right)+O(\mu^6).
\end{align}

So we have both $v$ (thanks to (\ref{eq:9})) and $p$ (thanks to
(\ref{eq:10})) and may rewrite (\ref{eq:8})$_1$ with the only fields
$u$ and $\eta$:
\begin{multline}
\label{eq:11}
u_t+u u_x+u_y \left(v|_{\varepsilon \eta}-\int_{\varepsilon \eta}^y u_x \right)+ \varepsilon \eta_x-\mu^2 \Bigg[(y-\varepsilon \eta)\left( (v|_{\varepsilon\eta})_t+u_x|_{\varepsilon\eta} \varepsilon \eta_t \right)\Bigg.\\
+\left(\int_{\varepsilon \eta}^y u \right) \left( (v|_{\varepsilon\eta})_x+u_x|_{\varepsilon\eta} \varepsilon \eta_x \right)-\left(\int_{\varepsilon \eta}^y u_x \right) v|_{\varepsilon \eta} -\int_{\varepsilon \eta}^y\int_{\varepsilon \eta}^{y'}u_{xt}\\
\left. -\int_{\varepsilon \eta}^y\left( u\int_{\varepsilon \eta}^{y'}u_{xx} \right) +\int_{\varepsilon \eta}^y\left(u_x \int_{\varepsilon \eta}^{y'} u_x \right)\right]_x= O\left(\frac{(u_y|_{\varepsilon \eta})_x}{\Rey} \right)+O\left( \frac{\mu}{\Rey} \right)+O\left(\frac{u_{yy}}{\mu \Rey}\right).
\end{multline}

In order to take off the dependence on $y$ of this equation,
we integrate between $y=-1+\mu^2 \gamma_{\infty}$ and $y=\varepsilon
\eta(x,t)$ and we define:
\begin{equation}
\label{eq:12}
\Hronde_{\mu,\gamma_{\infty}}=  1+\varepsilon \eta -\mu^2 \gamma_{\infty}, \mbox{ and }
\bar{u}(x,t)= \frac{1}{\Hronde_{\mu,\gamma_{\infty}}}\int_{-1+\mu^2 \gamma_{\infty}}^{\varepsilon \eta} u(x,y)\, {\rm d}y.
\end{equation}%

We also need a lemma that will enable to commute the integration and the $x$ differentiation under an assumption:

\begin{lemma}
\label{lemma1}
Let $F$ a $\mathcal{C}^1$ function defined in $\Omega =\{ (x,y) / x\in
\mathbb{R}, -1+\mu^2 \gamma_{\infty} < y < \varepsilon \eta(x) \} $, such that if $\forall x, \; F(x, y=\varepsilon \eta)=0,$ then
\begin{equation}
\label{eq:12.5}
\int_{-1+\mu^2 \gamma_{\infty}}^{\varepsilon \eta} \pd{F}{x}(x,y) {\rm d}y = \pd{ }{x} \int_{-1+\mu^2 \gamma_{\infty}}^{\varepsilon \eta} F(x,y) {\rm d} y.
\end{equation}%
\end{lemma}

The proof is very simple and left to the interested reader.

Thanks to Lemma \ref{lemma1}, one may commute the $x$ differentiation
of the square bracket in Equation (\ref{eq:11}) with the integral
since the terms in the square brackets vanish at $y=\varepsilon \eta$.
An integration by parts of the $\int u_y \left( v|_{\varepsilon
  \eta}-\int_{\varepsilon \eta}^y u_x\right) {\rm d}y$ term, and the
treatment of $\int (u^2)_x$ leads to (below, we write
$\Hronde=\Hronde_{\mu,\gamma_{\infty}}$):
%
\begin{multline}
\label{eq:14}
\Hronde \bar{u}_t+\left( \int_{-1+\mu^2 \gamma_{\infty}}^{\varepsilon \eta} u^2 \right)_x+\Hronde \varepsilon \eta_x+(\bar{u}-u|_{-1+\mu^2 \gamma_{\infty}})\left( \varepsilon \eta_t+\left( \Hronde \bar{u} \right)_x \right)-\bar{u}\left( \Hronde \bar{u} \right)_x\\
-\mu^2 \left[-\frac{\Hronde^2}{2}\left( (\partial_t+\bar{u} \partial_x)\left(v|_{\varepsilon \eta}\right)+\varepsilon \eta_t(u_x|_{\varepsilon \eta} -\bar{u}_x) +\varepsilon \eta_x(\bar{u}u_x|_{\varepsilon \eta}-\bar{u}_xu|_{\varepsilon \eta})\right) \right.\\
+\int_{-1+\mu^2 \gamma_{\infty}}^{\varepsilon \eta} \int_{\varepsilon \eta}^y (u-\bar{u}){\rm d} y'{\rm d}y \times \left( (v|_{\varepsilon \eta})_x+u_x|_{\varepsilon \eta} \varepsilon \eta_x \right)-\int_{-1+\mu^2 \gamma_{\infty}}^{\varepsilon \eta}\int_{\varepsilon \eta}^y (u-\bar{u})_x {\rm d}y'\, {\rm d}y \, v|_{\varepsilon \eta} \\
\left.-\int_{-1+\mu^2 \gamma_{\infty}}^{\varepsilon \eta} \int_{\varepsilon \eta}^y \left[ \int_{\varepsilon \eta}^{y'} u_{xt} + u\int_{\varepsilon \eta}^{y'} u_{xx}-u_x \int_{\varepsilon \eta}^{y'} u_x \right]{\rm d}y' {\rm d}y \right]_x= O\left(\frac{(u_y|_{\varepsilon \eta})_x}{\Rey} \right)+O\left(\frac{u_{yy}}{\mu \Rey}\right).
\end{multline}

We need now the following (double) assumption:
\begin{equation}
  \label{eq:13.bis}
u(x,y,t)=\bar{u}(x,t)+\mu^2 \tilde{u}(x,y,t),
\end{equation}%
with
\begin{equation}
  \label{eq:13.ter}
\int_{-1+\mu^2 \gamma_{\infty}}^{\varepsilon \eta} \tilde{u}=0 \mbox{ and } \bar{u}(x,t)=\frac{1}{\Hronde_{\mu,\gamma_{\infty}}}\int_{-1+\mu^2 \gamma_{\infty}}^{\varepsilon \eta} u(x,y)\, {\rm d}y.
\end{equation}
The mean $\bar{u}$ is the same as before. Notice that the expansion of
a function around its mean value $\bar{u}$ is not an assumption. The
real assumption is that the discrepancy with the mean is small
($O(\mu^2)$). An other way to formulate this assumption is to look at
an expansion in $\mu^2$, in which one assumes that the zeroth order
term does not depend on $y$ and that the next order term is
zero-mean value. This gives two different ways to
see its consequences. Last but not least, this assumption is proved to
be true in Lemma~11 (Eq.~(77)) of \cite{LeMeur2015} in case of a Boussinesq
flow (where $\varepsilon$ is small) without the assumption of
irrotationality in the Euler part of the flow. We remind the reader
that we still assume we solve the Euler equations and not yet the
Navier-Stokes ones. So we are coherent.

Upon this assumption, (\ref{eq:14}) simplifies to:
\begin{multline}
\label{eq:15}
  \Hronde \bar{u}_t+\Hronde \bar{u} \bar{u}_x+\Hronde \Hronde_x+(\bar{u}-u|_{-1+\mu^2 \gamma_{\infty}})\left( \Hronde_t+\left( \Hronde \bar{u} \right)_x \right)\\
  -\mu^2 \left[-\frac{\Hronde^2}{2}( \partial_t+\bar{u} \partial_x)\left(v|_{\varepsilon \eta}\right) -\frac{\Hronde^3}{6}(\bar{u}_{xt}+\bar{u}\bar{u}_{xx}-\bar{u} \bar{u}_x) \right]_x  = O\left(\frac{(u_y|_{\varepsilon \eta})_x}{\Rey} \right)+O\left(\frac{u_{yy}}{\mu \Rey}\right) + O(\mu^4).
\end{multline}

\begin{remark}\label{rem3}
The attention may be drawn to the fact that
\[
\Hronde_t +\left(\Hronde \bar{u}\right)_x   = \varepsilon \eta_t + \Hronde_x \bar{u} +\Hronde \bar{u}_x = v|_{\varepsilon \eta}+\Hronde \bar{u}_x + O(\mu^2) = v|_{-1+\mu^2 \gamma_{\infty}} +O(\mu^2).
\]
In the Euler case, $v|_{-1+\mu^2 \gamma_{\infty}}=0$ since the flow does not cross the boundary. So we would not need to compute $u|_{-1+\mu^2 \gamma_{\infty}}$.
\end{remark}

\subsection{Resolution in the boundary layer}
\label{sec.BL}

We write the system that applies in the layer, extracted from (\ref{eq:NS_adim}):
\begin{equation}
\label{eq:16}
\left\{
\begin{aligned}
u_t+u u_x+ vu_y-\mu u_{xx}/\Rey- u_{yy}/(\mu \Rey)+p_x&=0,\\
\mu^2\left(v_t+u \, v_x+v \, v_y\right)-\mu^3 v_{xx}/\Rey-\mu  v_{yy}/\Rey+p_y&=0,\\
u_x+v_y &= 0,\\
u(y=-1)=v(y=-1) &=0.
\end{aligned}
\right.
\end{equation}

This system may be rewritten with the change of variables justified in
(\ref{eq6.55}) $ y=-1+\mu^2 \, \gamma$, where $\gamma >0$ and may be
up to a large (but not too large) $\gamma_{\infty}$. We also use the
assumption (\ref{assumption_Re}) on $\Rey$ such that
$\Rey= R \, \mu^{-5}$ where $R$ is a constant. We should have tilded
the fields but would have dropped the tilde soon after. So we omit
them. When precision is needed, we denote $u^{BL}=u(x,\gamma,t)$ the
horizontal velocity in the boundary layer. The system writes:
\begin{equation}\label{eq:16.bis}
\left\{
\begin{aligned}
u_t+u\, u_x+v\, u_{\gamma}/\mu^2 -\mu^6/R\, u_{xx}-u_{\gamma \gamma}/R+p_x &=0,\\
\mu^2 \left(v_t+u\, v_x +v\, v_{\gamma}/\mu^2 \right)-(\mu^8/R) v_{xx}-(\mu^2/R) v_{\gamma \gamma}+p_{\gamma}/\mu^2 & = 0,\\
u_x+v_{\gamma}/\mu^2 &= 0, \\
u(x,\gamma=0,t)=v(x,\gamma=0,t)&=0.
\end{aligned}
\right.
\end{equation}

As is classical, we first compute $v$ (owing to (\ref{eq:16.bis}$_{3}$) and (\ref{eq:16.bis}$_{4}$):
\begin{equation}
\label{eq:17}
v(x,\gamma,t)=0-\mu^2 \int_{\gamma=0}^{\gamma } u_x \, {\rm d}\gamma'.
\end{equation}
Then we can compute the differentiated pressure from
(\ref{eq:16.bis})$_2$ that proves $p_{\gamma}=O(\mu^4)$. As a
consequence,
\[
p^{BL}(x,\gamma,t)=p^{BL}(\gamma \rightarrow \gamma_{\infty})+O(\mu^4),
\]
where $p^{BL}(\gamma \rightarrow \gamma_{\infty})$ is determined
thanks to a matching condition with the bottom of the upper part
(Euler part). From (\ref{eq:10}), and owing to the already stated
assumption (\ref{eq:13.bis}, \ref{eq:13.ter}) the pressure in the
boundary layer is, up to $O(\mu^4)$:
\begin{equation}
\label{eq:20}
p^{BL}(x,\gamma,t)=\varepsilon \eta(x,t)+\mu^2\left[\Hronde(\partial_t+\bar{u}\partial_x)(v|_{\varepsilon \eta})+\Hronde^2/2(\bar{u}_{xt}+\bar{u} \bar{u}_{xx}-\bar{u}_x\bar{u}_x)\right]+O(\mu^4).
\end{equation}
Last, we may gather $v^{BL}$ (from (\ref{eq:17})), $p^{BL}$ (from
(\ref{eq:20})) and rewrite (\ref{eq:16.bis})$_1$:
\begin{multline}\label{eq:21}
  u^{BL}_t+u^{BL}\, u^{BL}_x-u^{BL}_{\gamma}\int_0^{\gamma}u^{BL}_x(\gamma'){\rm d}\gamma'-\frac{u^{BL}_{\gamma \, \gamma}}{R} +\varepsilon \eta_x\\
  +\mu^2\left[\Hronde (\partial_t+\bar{u}\partial_x)(v|_{\varepsilon \eta})+\Hronde^2/2(\bar{u}_{xt}+\bar{u} \bar{u}_{xx}-\bar{u}_x\bar{u}_x)\right]_x=O(\mu^4).
\end{multline}
At that stage of the derivation, we recognize a Prandtl's equation. We
do not know how to derive a simpler model. It is well-known that
Prandtl's equation still resists to the best of physicists and
mathematicians. It was proved to be ill-posed in
\cite{Gerard-Varet_Dormy_2010} and partially well-posed later. 
Moreover, so as to couple this
equation with (\ref{eq:15}), one should assume a link between $u^{BL}$
and $u|_{-1+\mu^2 \gamma_{\infty}}$ like identity by continuity.

\section{Conclusions}

We stopped our derivation, in the fully nonlinear regime, at equation
(\ref{eq:15}), (\ref{eq:21}), which is only an intermediate state
because we still have functions of $x,\gamma,t$. We would like to stress it out that
Equations (\ref{eq:15}), (\ref{eq:21}) are still Galilean invariant despite the presence
of the boundary layer. The proposed model enjoys this property because we did not
introduce any drastic simplifications yet at this level. To make further
progress, the Prantdl-type equation should be further simplified but it
seems highly speculative. One strategy could consist in assuming a
particular profile of the velocity $u^{BL}$ in the coordinate
$\gamma$ similar to the one computed in \cite{LeMeur2015} in the
Boussinesq regime, but it is incoherent. Further research is needed
to reach an effective 1D model.

\vspace*{-1.0em}
\section*{Acknowledgements}

This research did not receive any specific grant from funding agencies in the public, commercial, or not-for-profit sectors.


\end{document}